\documentclass[letterpaper, 10 pt, conference]{ieeeconf}

\pdfoutput=1

\IEEEoverridecommandlockouts
\usepackage{amsmath} 
\usepackage{amssymb}  
\usepackage[utf8]{inputenc}
\usepackage{algorithm2e}
\usepackage{color}
\usepackage{tikz}
\usetikzlibrary{shapes,arrows.meta}
\usepackage{subfig}

\usepackage{graphicx}
\usepackage[font=footnotesize,labelfont=bf]{caption}

\DeclareMathOperator*{\argmin}{arg\,min}
\DeclareMathOperator{\prox}{prox}
\DeclareMathOperator{\refl}{refl}

\newtheorem{defn}{Definition}
\newtheorem{remark}{Remark}
\newtheorem{prop}{Proposition}

\newtheorem{assumption}{Assumption}

\def\oprocend{\hfill$\square$}
\def\oprocendbis{\hfill$\blacksquare$}

\title{\LARGE \bf
Distributed Optimization over Lossy Networks via Relaxed Peaceman-Rachford Splitting: a Robust ADMM Approach
}
\author{N. Bastianello$^\dagger$, M. Todescato$^\ddagger$, R. Carli$^\dagger$, L. Schenato$^\dagger$
\thanks{$^\dagger$ Department of Information Engineering (DEI), University of Padova, Italy. {\tt\small nicola.bastianello.3@studenti.unipd.it, [carlirug|schenato]@dei.unipd.it}.}
\thanks{$^\ddagger$ Bosch Center for Artificial Intelligence. Renningen, Germany. {\tt\small mrc.todescato@gmail.com}. The work was carried out during the author's postdoctoral fellowship at DEI.}
}

\begin{document}

\maketitle
\thispagestyle{empty}
\pagestyle{empty}

\begin{abstract}
In this work we address the problem of distributed optimization of the sum of convex cost functions in the context of multi-agent systems over lossy communication networks. Building upon operator theory, first, we derive an ADMM-like algorithm that we refer to as relaxed ADMM (R-ADMM) via a generalized \emph{Peaceman-Rachford Splitting} operator on the Lagrange dual formulation of the original optimization problem. This specific algorithm depends on two parameters, namely the averaging coefficient $\alpha$ and the augmented Lagrangian coefficient $\rho$. We show that by setting $\alpha=1/2$ we recover the standard ADMM algorithm as a special case of our algorithm. Moreover, by properly manipulating the proposed R-ADMM, we are able to provide two alternative ADMM-like algorithms that present easier implementation and reduced complexity in terms of memory, communication and computational requirements. Most importantly the latter of these two algorithms provides the first ADMM-like algorithm which has guaranteed convergence even in the presence of lossy communication under the same assumption of standard ADMM with lossless communication. Finally, this work is complemented with a set of compelling numerical simulations of the proposed algorithms over cycle graphs and random geometric graphs subject to i.i.d. random packet losses.
\end{abstract}


\section{Introduction}\label{sec:intro}
From classical control theory to more recent Machine Learning applications, many problems can be cast as optimization problems \cite{Slavakis:14} and, in particular, as large-scale optimization problems given the advent of Internet-of-Things we are witnessing with its ever-increasing growth of large-scale cyber-physical systems. Hence, stemming from classical optimization theory, in order to break down the computational complexity, parallel and distributed optimization methods have been the focus of a wide branch of research \cite{BD:TP:1989}. Within this vast topic, typical applications, going under the name of \emph{distributed consensus optimization}, foresee distributed computing nodes to communicate in order to achieve a desired common goal. More formally, the distributed nodes seek to
$$
\min_{x}\sum_{i=1}^Nf_i(x)
$$
where, usually, each $f_i$ is owned by one node only.
Toward this application among very many different optimization algorithms explored in past as well as in current literature, e.g. subgradient methods \cite{BJ:MR:MJ:2010}, the well known Alternating Direction Method of Multipliers (ADMM), 
%
first introduced in \cite{glowinski1975approximation} and \cite{gabay1976dual}, is recently receiving an ever-increasing interest because of its numerical efficiency and its natural structure which makes it well-suited for distributed and parallel computing. In particular, the relatively recent monograph \cite{boyd2011distributed} reveals the ADMM in detail presenting a broad set of selected applications to which ADMM is suitably applied. For a wider set of applications together with some convergence results we refer the interested reader to \cite{fukushima1992application,eckstein1994some,eckstein1992douglas,EG:AT:ES:MJ:2015}.\\
%
While ADMM can be proficiently applied to distributed setups, rigorous convergence results are usually provided only in scenarios characterized by synchronous updates and lossless communications.However, practically it is rarely possible and often difficult to ensure synchronization and communication reliability among computing nodes. And even when this is possible via specific communication protocols, it is clear how the impossibility to deal with asynchronous and lossy updates would majorly compromise the algorithm applicability.\\
Hence, an extensive body of work has been devoted to overcoming this limitation by adapting the ADMM to operate in an asynchronous fashion. Among the first steps in this direction, \cite{iutzeler2013asynchronous} proves convergence when only one randomly selected coordinate is updated at each iteration. Similarly, \cite{wei20131} suggests to update only the variables related to a subset of constraints randomly selected at each iteration, showing convergence of the algorithm with a rate of $O(1/k)$. To deal with asynchronous updates, a master-slave architecture is proposed in \cite{zhang2014asynchronous}.
The more recent \cite{bianchi2016coordinate} extends the formulation introduced in \cite{iutzeler2013asynchronous} to allow the update of a subset of coordinates at each instant. In \cite{chang2016asynchronous}, in view of large-scale optimization, the convergence rate of a partially asynchronous ADMM -- i.e., subject to a maximum allowed delay  -- is studied. Finally \cite{peng2016arock} defines a framework for asynchronous operations used to solve a broad class of optimization problems and showing how to derive an asynchronous ADMM formulation.\\
%
Conversely to the above works that deal with asynchronous updates, to the best of our knowledge, no work explicitly focuses on the robustness of ADMM to packet losses. Yet, to set the stage for the analysis of robustness of the ADMM algorithm to losses in the communication, we resort to a different body of literature on operator theory. Here, the underlying idea is to convert optimization problems into the problem of finding the fixed points of suitable nonexpansive operators \cite{bauschke2011convex}. However, the mere application of the so-called \emph{proximal point algorithm} (PPA)  -- introduced in \cite{rockafellar1976monotone} and the later \cite{parikh2014proximal} -- to look for the fixed points can be unwieldy in complex optimization problems. Hence, particular credits have been given to \emph{splitting methods} which exploit the problem's structure to break it in smaller and more manageable pieces. It is in the framework of splitting operators, and in particular the well recognized Peaceman-Rachford (PRS) \cite{peaceman1955numerical} and Douglas-Rachford (DRS) \cite{douglas1956numerical,lions1979splitting} splitting, that the ADMM comes into place. Indeed, the classical formulation of the ADMM naturally arises as application of the DRS to the Lagrange dual problem of the original optimization problem \cite{eckstein2012augmented}. For further details on a variety of splitting operators and their application in asynchronous setups we refer to \cite{davis2016convergence} and \cite{hannah2016unbounded}, respectively.\\
%
In this paper we present and analyze different formulation for the ADMM algorithm. We are particularly interested to the broad class of distributed consensus optimization problems. Our final goal is to present a novel robustness result in scenarios characterized by synchronous but possibly lossy updates among distributed nodes. To achieve our result we start by considering a prototypical optimization problem assuming reliable loss-free communication. In this case, by leveraging the general framework arising from the Krasnosel'skii-Mann (KM) iteration for averaging operators \cite{krasnosel1955two,mann1953mean}, we derive a relaxed version of the ADMM (R-ADMM). Next, we draw the attention to the problem of interest, i.e., distributed consensus optimization. We first present the natural algorithmic implementation of the R-ADMM tailored for the problem. Then we propose two different implementations which are particularly favorable for storage and communication purposes. Moreover, the latter turns out extremely advantageous and yet robust in the presence of lossy communication.
As natural byproduct we obtain a comprehensive and self-contained overview on the algorithm and a plethora of possible practical implementations.\\
%
The remainder of the paper is organized as follows. Section~\ref{sec:operators-background} presents the necessary background on splitting operators. Section~\ref{sec:ADMMandRADMM} reviews the classical ADMM algorithm and its generalized version. Section~\ref{sec:distributed_consensus} focuses on the analysis of distributed consensus optimization. Section~\ref{sec:simulation} collects some numerical simulations and Section~\ref{sec:conclusions} concludes the paper. The technical proofs can be found in the Appendices.

\section{Background on Splitting Operators}\label{sec:operators-background}

This Section introduces some background on operator theory on Hilbert spaces and, in particular, on nonexpansive operators. The interest for operator theory stems from the fact that a convex optimization problem can be cast into the problem of finding the fixed point(s) of a suitable nonexpansive operator $T$ \cite{davis2016convergence,bauschke2011convex}, that is the points $x^*$ such that $Tx^*=x^*$. 

\subsection{Definitions and Properties}\label{subsec:definitions}
\begin{defn}[Nonexpansive operators]
Let $\mathcal{X}$ be a Hilbert space, an operator $T:\mathcal{X}\rightarrow\mathcal{X}$ is said to be \textit{nonexpansive} if it has unitary Lipschitz constant, \textit{i.e.} it verifies $\|Tx-Ty\|\leq\|x-y\|$ for any two $x,y\in\mathcal{X}$.\oprocend
\end{defn}

\smallskip


\begin{defn}[$\alpha$-averaged operators]
Let $\mathcal{X}$ be a Hilbert space, $T:\mathcal{X}\rightarrow\mathcal{X}$ a nonexpansive operator and $\alpha\in(0,1)$. We define the \textit{$\alpha$-averaged operator} $T_\alpha$ as $T_\alpha=(1-\alpha)I+\alpha T$, where $I$ is the identity operator on $\mathcal{X}$.\oprocend
\end{defn}

\smallskip

Notice that $\alpha$-averaged operators are also nonexpansive, indeed nonexpansive operators are $1$-averaged. Moreover, the $\alpha$-averaged operator $T_\alpha$ has the same fixed points of $T$ \cite{bauschke2011convex}.


\smallskip

\begin{defn}[Proximal and reflective operators]
Let $\mathcal{X}$ be a Hilbert space and $f:\mathcal{X}\rightarrow\mathbb{R}\cup\{+\infty\}$ be a closed, proper and convex function. We define the \textit{proximal operator} of $f$ with penalty $\rho>0$, $\prox_{\rho f}:\mathcal{X}\rightarrow\mathcal{X}$, as
\begin{equation*}
	\prox_{\rho f}(y)=\argmin_{x\in\mathcal{X}}\left\{f(x)+\frac{1}{2\rho}\|x-y\|^2\right\}.
\end{equation*}
Moreover, we define the relative \textit{reflective operator} as $\refl_{\rho f}=2\prox_{\rho f}-I$.\oprocend
\end{defn}

\smallskip

It can be seen that the proximal operator is $1/2$-averaged and the reflective operator is nonexpansive \cite{bauschke2011convex}.

\subsection{Finding the Fixed Points of Nonexpansive Operators}\label{subsec:fixpoints}
One of the prototypical algorithm for finding the fixed points of $T$ is the Krasnosel'skii-Mann (KM) iteration \cite{bauschke2011convex}
\begin{equation}\label{eq:km-iteration}
	x(k+1)=T_\alpha x(k)=(1-\alpha)x(k)+\alpha Tx(k)
\end{equation}
where in general the step size $\alpha$ can be time-varying. 
Notice that the KM iteration is equivalent to $x(k+1)=x(k)-\alpha Sx(k)$, where $S=I-T$, that is, finding the fixed points of $T$ coincides with finding the zeros of $S$.\\
Now, consider the general unconstrained problem
\begin{equation}\label{eq:MinProblem}
\min_{x\in\mathcal{X}} \{f(x)+g(x)\},
\end{equation} 
where $f,g$ are closed proper and convex not necessarily smooth functions. Further, assume that simultaneous minimization of $f+g$ is unwieldy while minimizing $f$ and $g$ separately is manageable. In this case, to compute the solution of \eqref{eq:MinProblem}
we can apply the KM iteration to the \textit{Peaceman-Rachford Splitting operator}, defined as (see \cite{davis2016convergence,peng2016arock}),
\begin{equation}
	T_{PRS}=\refl_{\rho f}\circ\refl_{\rho g}.
\end{equation}
As show in \cite{bauschke2011convex}, the iteration 
\begin{equation}\label{eq:KM_T_PRS}
	x(k+1)=(1-\alpha)x(k)+\alpha T_{PRS}x(k)
\end{equation}
can be conveniently implemented by the following updates
\begin{align}
	\psi(k)&=\prox_{\rho g}(z(k))\label{eq:prs-1}\\
	\xi(k)&=\prox_{\rho f}(2\psi(k)-z(k))\label{eq:prs-2}\\
	z(k+1)&=z(k)+2\alpha(\xi(k)-\psi(k))\label{eq:prs-3}
\end{align}
where $\psi, \xi, z$ are suitable auxiliary variables while the optimal solution $x^*$ to \eqref{eq:MinProblem} is recovered from the limit $z^*$ of the iterate $z(k)$ by computing $x^*=\prox_{\rho g}(z^*)$. This algorithm goes under the name of \textit{relaxed Peaceman-Rachford splitting} (R-PRS), where ``relaxed'' denotes the fact that the KM iteration is $\alpha$-averaged. In case $\alpha=1$ we recover the classic \textit{Peaceman-Rachford splitting} introduced in \cite{peaceman1955numerical}, and in case $\alpha=1/2$ we recover the \textit{Douglas-Rachford splitting} \cite{douglas1956numerical}.\\
The important feature of splitting schemes such as the R-PRS is that they divide the computational load of iterate \eqref{eq:KM_T_PRS} into smaller subproblems that can be solved more efficiently.

\section{From the ADMM to the Relaxed ADMM}\label{sec:ADMMandRADMM}
In this Section, we first review the popular ADMM algorithm \cite{gabay1976dual,boyd2011distributed}, then we introduce the more general \emph{relaxed} ADMM (R-ADMM) algorithm, and compare the two methods.

\subsection{The ADMM Algorithm}\label{subsec:ADMM}
Consider the following optimization problem
\begin{align}\label{eq:primal-problem}
\begin{split}
	&\min_{x\in\mathcal{X},y\in\mathcal{Y}} \{f(x)+g(y)\}\\
	&\text{s.t.}\ Ax+By=c
\end{split}
\end{align}
where $\mathcal{X}$ and $\mathcal{Y}$ are Hilbert spaces, $f:\mathcal{X}\rightarrow\mathbb{R}\cup\{+\infty\}$ and $g:\mathcal{Y}\rightarrow\mathbb{R}\cup\{+\infty\}$ are closed, proper and convex functions\footnote{A function $f:\mathcal{X}\rightarrow\mathbb{R}\cup\{+\infty\}$ is said to be \textit{closed} if $\forall a \in\mathbb{R}$ the set $\{x\in\operatorname{dom}(f)\ |\ f(x)\leq a\}$ is closed. Moreover, $f$ is said to be \textit{proper} if it does not attain $-\infty$ \cite{boyd2011distributed}.}.\\
To solve problem \eqref{eq:primal-problem} via the ADMM algorithm, we first define the \textit{augmented Lagrangian} as
\begin{align}\label{eq:augmented-lagr}
\begin{split}
	\mathcal{L}_{\rho}(x,y;w)&=f(x)+g(y)-w^\top\left(Ax+By-c\right)\\&+\frac{\rho}{2}\|Ax+By-c\|^2
\end{split}
\end{align}
where $\rho>0$ and $w$ is the vector of Lagrange multipliers. The ADMM algorithm consists in keeping alternating the following update equations
\begin{align}
	y(k+1)&=\operatorname{arg\,min}_y \mathcal{L}_\rho(x(k),y;w(k))\label{eq:admm-1}\\
	w(k+1)&=w(k)-\rho(Ax(k)+By(k+1)-c)\label{eq:admm-2}\\
	x(k+1)&=\operatorname{arg\,min}_x \mathcal{L}_\rho(x,y(k+1);w(k+1))\label{eq:admm-3}.
\end{align}
Notice that the above formulation is equivalent to the one proposed in \cite{boyd2011distributed} except for a change in the order of the updates which however does not affect the convergence properties of the algorithm. Moreover, the ADMM algorithm is provably shown to converge to the optimal solution of \eqref{eq:primal-problem} for any $\rho>0$ assuming that $\mathcal{L}_0$ has a saddle point \cite{boyd2011distributed}.\\
We conclude this section by stressing the following fact. While the ADMM in its classical form \eqref{eq:admm-1}--\eqref{eq:admm-3} is typically presented as an augmented Lagrangian method computed with respect to the primal problem~\eqref{eq:primal-problem}, the algorithm naturally arises from the application of the DRS to the Lagrange dual of problem~\eqref{eq:primal-problem}. This will be made clear in the next section.

\subsection{The Relaxed ADMM}\label{subsec:R-ADMM}
The Relaxed ADMM algorithm can be derived applying the R-PRS method described in Section \ref{sec:operators-background} to the Lagrange dual of problem \eqref{eq:primal-problem}, that is to
\begin{equation}\label{eq:dual-problem}
	\min_{w\in\mathcal{W}}\left\{d_f(w)+d_g(w)\right\}
\end{equation}
where
\begin{align*}
	d_f(w)&=f^*(A^\top w)\\
	d_g(w)&=g^*(B^\top w)-w^\top c,
\end{align*}
and $f^*$, $g^*$ are the convex conjugates of $f$ and $g$\footnote{The \textit{convex conjugate} of a function $f$ is defined as $f^*(y)=\sup_{x\in\mathcal{X}}\{\langle y,x\rangle-f(x)\}$.}. The derivation of problem \eqref{eq:dual-problem} can be found in \cite{davis2016convergence,peng2016arock}.\\ 
Observe that, given the structure of problem \eqref{eq:primal-problem}  (i.e., proper closed and convex functions and linear constraints) there is no duality gap and, in turn, the optimal solutions of \eqref{eq:primal-problem} and of \eqref{eq:dual-problem} attain the same optimal value.\\
The motivation for dealing with the Lagrange dual problem relies on the fact that the minimization in \eqref{eq:dual-problem} is performed over a single variable, thus allowing for the use of the R-PRS algorithm described in \eqref{eq:prs-1}, \eqref{eq:prs-2} and \eqref{eq:prs-3}.\\
Lemma $11$ in \cite{davis2016convergence} shows that the update \eqref{eq:prs-1} and the update \eqref{eq:prs-2}, applied to the dual problem, can be conveniently computed by, respectively,
\begin{align}
	y(k)&=\argmin_y\left\{g(y)-z^\top(k)(By-c)+\frac{\rho}{2}\|By-c\|^2\right\}\nonumber\\
	\psi(k)&=z(k)-\rho(By(k)-c)\label{eq:psi-update}
	\end{align}
	and
	\begin{align}
	x(k)&=\argmin_x\left\{f(x)-(2\psi(k)-z(k))^\top Ax+\frac{\rho}{2}\|Ax\|^2\right\}\nonumber\\
	\xi(k)&=2\psi(k)-z(k)-\rho Ax(k)\label{eq:xi-update}
\end{align}
The so called \emph{relaxed ADMM} algorithm (in short R-ADMM) consists in applying iteratively the set of five equations given by the two equations in \eqref{eq:psi-update}, the two equations in \eqref{eq:xi-update} and equation \eqref{eq:prs-3}.\\
It is worth stressing a fundamental difference regarding the auxiliary variables $z$ in \eqref{eq:prs-1}--\eqref{eq:prs-3} and those used in \eqref{eq:psi-update},\eqref{eq:xi-update}. Indeed, when implementing the R-PRS \eqref{eq:prs-1}--\eqref{eq:prs-3}, the KM iteration is applied directly to the primal problem \eqref{eq:MinProblem}. Hence $z$ has the same dimension of the primal variable $x$. Conversely, the R-ADMM as in \eqref{eq:psi-update}, \eqref{eq:xi-update} and \eqref{eq:prs-3}, is derived on the dual problem \eqref{eq:dual-problem}. Hence, in this case $z$ has dimension of the constraints. This will become more clear later when dealing with distributed consensus optimization problems.\\
%
%
Next, we derive a  more compact formulation of the R-ADMM, that shows clearly the relation with the popular ADMM algorithm described in \eqref{eq:admm-1}, \eqref{eq:admm-2} and \eqref{eq:admm-3}.\\
First of all, by adding $\psi(k)$ on both sides of the second equation in \eqref{eq:psi-update},
we obtain
\begin{equation}
	2\psi(k)-z(k)=\psi(k)-\rho(By(k)-c)\label{eq:equality-1}
\end{equation}
and substituting this equation in \eqref{eq:xi-update} we get
\begin{equation}
	\xi(k)=\psi(k)-\rho(Ax(k)+By(k)-c)\label{eq:equality-2}.
\end{equation}
Getting $z(k)$ from \eqref{eq:equality-1} and, substituting back into the first equation in \eqref{eq:xi-update} we obtain
\begin{align*}
\begin{split}
	x(k)&=\argmin_x\{f(x)-\psi^\top(k)(Ax+By(k)-c)\\&+\frac{\rho}{2}\|Ax+By(k)-c\|^2\}.
\end{split}
\end{align*}
where terms independent of $x$ were added.
By substituting \eqref{eq:equality-1} and \eqref{eq:equality-2} in \eqref{eq:prs-3} we get
\begin{equation}
	z(k+1)=\psi(k)-\rho Ax(k)-\rho(2\alpha-1)(Ax(k)+By(k)-c)\label{eq:z-update}
\end{equation}
and by plugging \eqref{eq:z-update} into the first equation in \eqref{eq:psi-update} and adding some terms that do not depend on $y$, we get
\begin{align*}
\begin{split}
	y&(k+1)=\argmin_y\{g(y)\\
	&-\psi^\top(k)(Ax(k)+By-c)+\rho\|By-c\|^2\\
	&+\rho\left[Ax(k)+(2\alpha-1)(Ax(k)+By(k)-c)\right]^\top(By-c)\}.
\end{split}
\end{align*}
Finally, recalling the defintion of augmented Lagrangian \eqref{eq:augmented-lagr} and renaming $\psi$ as $w$, we arrive to the three updates that represent the R-ADMM algorithm
\begin{align}
\begin{split}
	y(k+1)&=\argmin_y\{\mathcal{L}_\rho(x(k),y;w(k))\\&+\rho(2\alpha-1)\langle By,(Ax(k)+By(k)-c)\rangle\}\label{eq:r-admm-1}
\end{split}\\
\begin{split}
	w(k+1)&=w(k)-\rho(Ax(k)+By(k+1)-c)\\&-\rho(2\alpha-1)(Ax(k)+By(k)-c)\label{eq:r-admm-2}
\end{split}\\
	x(k+1)&=\argmin_x\mathcal{L}_\rho(x,y(k+1);w(k+1)).\label{eq:r-admm-3}
\end{align}
The ADMM algorithm in \eqref{eq:admm-1}--\eqref{eq:admm-3} can be recovered from this formulation of the R-ADMM by setting $\alpha=1/2$, which cancels the additional terms weighted by $2\alpha-1$.\\
It is of notice that the R-ADMM has two tunable parameters, $\rho$ and $\alpha$, against the only one of the ADMM, $\rho$, which is the cause of the greater reliability of the R-ADMM.

\smallskip

\begin{remark}
Figure \ref{fig:relationships-algorithms} depicts the relationships between the splitting operators derived from the relaxed Peaceman-Rachford, and the classic and relaxed ADMM.
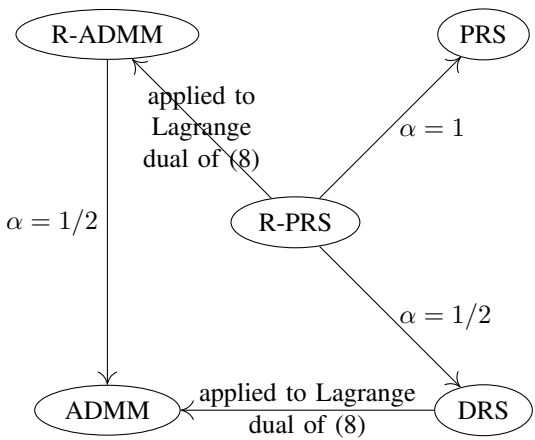
\begin{figure}[!ht]
\centering
\begin{tikzpicture}
	\node[draw,ellipse] (r-prs) at (0,0) {R-PRS};
	\node[draw,ellipse] (prs) at (2.5,2.5) {PRS};
	\node[draw,ellipse] (drs) at (2.5,-2.5) {DRS};
	\node[draw,ellipse] (r-admm) at (-2.5,2.5) {R-ADMM};
	\node[draw,ellipse] (admm) at (-2.5,-2.5) {ADMM};
	
	\path[-{>[scale=2.5,length=2,width=3]}]
	(r-prs) edge node[right] {$\alpha=1/2$} (drs)
	(r-prs) edge node[right] {$\alpha=1$} (prs)
	(r-prs) edge node[align=center] {applied to\\ Lagrange\\ dual of \eqref{eq:primal-problem}} (r-admm)
	(drs) edge node[align=center] {applied to Lagrange\\ dual of \eqref{eq:primal-problem}} (admm)
	(r-admm) edge node[left] {$\alpha=1/2$} (admm);
\end{tikzpicture}
\caption{Relationships between the algorithms.}
\label{fig:relationships-algorithms}
\end{figure}
\oprocend
\end{remark}

\section{Distributed Consensus Optimization}\label{sec:distributed_consensus}
This Section introduces the distributed consensus convex optimization problem that we are interested in, and the solutions obtained by applying the R-ADMM algorithm.

\subsection{Problem Formulation}\label{subsec:distributed_problem}
Let $\mathcal{G}=(\mathcal{V},\mathcal{E})$ be a graph, with $\mathcal{V}$ the set of $N$ vertices, labeled $1$ through $N$, and $\mathcal{E}$ the set of undirected edges. For $i \in \mathcal{V}$, by $\mathcal{N}_i$ we denote the set of neighbors of node $i$ in $\mathcal{G}$, namely,
$$
\mathcal{N}_i =\left\{j \in V \,:\, (i,j) \in \mathcal{E} \right\}.
$$
We are interested in solving the following optimization problem
\begin{align}\label{eq:opt_problem}
\begin{split}
	&\min_{x}\sum_{i=1}^Nf_i(x)
\end{split}
\end{align}
where $f_i:\mathbb{R}^n \rightarrow\mathbb{R}\cup\{+\infty\}$ are closed, proper and convex functions and where $f_i$ is known only to node $i$. In the following we denote by $x^*$ the optimal solution of \eqref{eq:opt_problem}.\\
Observe that \eqref{eq:opt_problem} can be equivalently formulated as
\begin{align}\label{eq:distributed-primal}
\begin{split}
	&\min_{x_i,\forall i}\sum_{i=1}^Nf_i(x_i)\\
	&\text{s.t.}\ x_i=x_j,\ \forall (i,j)\in\mathcal{E}
\end{split}
\end{align}
By introducing for each edge $(i,j)\in\mathcal{E}$ the two \textit{bridge variables} $y_{ij}$ and $y_{ji}$, the constraints in \eqref{eq:distributed-primal} can be rewritten as
\begin{align*}
\begin{split}
	& x_i=y_{ij}\\
	& x_j=y_{ji}\\
	& y_{ij}=y_{ji}
\end{split}\ \ \ \forall (i,j)\in\mathcal{E}.
\end{align*}
Defining $\mathbf{x}=[x_1^\top,\ldots,x_N^\top]^\top$, $f(\mathbf{x})=\sum_if_i(x_i)$, and stacking all bridge variables in $\mathbf{y} \in \mathbb{R}^{n |\mathcal{E}|}$, we can reformulate the problem as
\begin{align*}
	& \min_{\mathbf{x}} f(\mathbf{x})\\
	& \text{s.t.}\ \ A\mathbf{x}+\mathbf{y}=0\\
	& \mathbf{y}=P\mathbf{y}
\end{align*}
for a suitable $A$ matrix and with $P$ being a permutation matrix that swaps $y_{ij}$ with $y_{ji}$. Making use of the indicator function $\iota_{(I-P)}(\mathbf{y})$ which is equal to 0 if $(I-P)\mathbf{y}=0$, and $+\infty$ otherwise, we can finally rewrite problem \eqref{eq:distributed-primal} as
\begin{align}\label{eq:primal-indicator-f}
\begin{split}
	& \min_{\mathbf{x},\mathbf{y}}\left\{f(\mathbf{x})+\iota_{(I-P)}(\mathbf{y})\right\}\\
	& \text{s.t.}\ \ A\mathbf{x}+\mathbf{y}=0.
\end{split}
\end{align}
In next Section we apply the R-ADMM algorithm described in Section \ref{subsec:R-ADMM} to the above problem.
%


\subsection{R-ADMM for Convex Distributed Optimization}\label{subsec:distributed_R-ADMM}

In this section we employ \eqref{eq:r-admm-1}, \eqref{eq:r-admm-2} and \eqref{eq:r-admm-3} to solve problem \eqref{eq:primal-indicator-f}. 
To do so we introduce the dual variables $w_{ij}$ and $w_{ji}$ which are associated to the constraints $x_i=y_{ij}$ and $x_j=y_{ji}$, respectively.
The resulting algorithm is described in  Algorithm~\ref{alg:r-admm-three-eqs}. Observe that R-ADMM applied to \eqref{eq:primal-indicator-f} is amenable of a \emph{distributed} implementation, in the sense that node $i$ stores in memory only the variables $x_i$, $y_{ij}$, $w_{ij}, j\in\mathcal{N}_i$, and updates these variables exchanging information only with its neighbors, i.e, with nodes in $\mathcal{N}_i$. 
%

\begin{algorithm}[ht!]
	\SetKwInOut{Input}{Input}
	\Input{Set the termination condition $K>0$. For each node $i$, initialize $x_i(0)$, $\{y_{ij}(0),w_{ij}(0)\}_{j\in\mathcal{N}_i}$.}
	$k\leftarrow0$\;
	\While{$k<K$ every agent $i$}{
		($i$ collects $\{x_j(k),y_{ji}(k),w_{ji}(k)\}$ received from each neighbor $j\in\mathcal{N}_i$)\;
		compute in order 
		\begin{align*}
			&y_{ij}(k+1)=\frac{1}{2\rho}\Big[(w_{ij}(k)+w_{ji}(k))+\\ 
			&\ \ 2\alpha\rho(x_i(k)+x_j(k)) -\rho(2\alpha-1)(y_{ij}(k)+y_{ji}(k))\Big]
		\end{align*}	
		\begin{align*}
			&w_{ij}(k+1)=\frac{1}{2}\Big[(w_{ij}(k)-w_{ji}(k))+\\
			&\ \ 2\alpha\rho(x_i(k)-x_j(k)) -\rho(2\alpha-1)(y_{ij}(k)-y_{ji}(k))\Big]
		\end{align*}
		\begin{align*}
			x_i(k+1)&=\argmin_{x_i}\Bigg\{f_i(x_i)+\frac{\rho}{2}|\mathcal{N}_i| \|x_i\|^2 \\
			&+\Big(\sum_{j\in\mathcal{N}_i}w_{ij}(k+1)-\rho y_{ij}(k+1)\Big)^\top x_i\Bigg\};\;
		\end{align*}\\
		broadcast $x_i(k+1)$, $y_{ij}(k+1)$ and $w_{ij}(k+1)$ to each neighbor $j$\;
		$k\leftarrow k+1$\;
	}
\caption{Distributed R-ADMM using \eqref{eq:r-admm-1}--\eqref{eq:r-admm-3}.}
\label{alg:r-admm-three-eqs}
\end{algorithm}

Notice that, in the update of $x_i$, the term $w_{ij}(k+1)-\rho y_{ij}(k+1)$ can be rewritten, using the previous updates, as a function of the variables computed at time $k$ only. Therefore, only one round of transmissions is necessary.\\
The above implementation of the R-ADMM is quite straightforward and popular but very unwieldy due to the fact that, depending on the number of neighbors, there might be nodes which need to store, update and transmit a large number of variables. The derivation of Algorithm~\ref{alg:r-admm-three-eqs} is reported in Appendix~\ref{app:derivation-alg-3eqs}.

In the following we provide an alternative algorithm which is derived directly from the application of the set of five equations in \eqref{eq:psi-update}, \eqref{eq:xi-update} and \eqref{eq:prs-3} to the dual of problem \eqref{eq:primal-indicator-f}. Notice, that since the vector $z$ has the same dimension of the vector $w$, this implies the presence of also the variables $z_{ij}$ and $z_{ji}$ for any $(i,j) \in \mathcal{E}$.\\
We have the following Proposition, which is proved in Appendix~\ref{app:proof-prop-1}.

\smallskip

\begin{prop}\label{pr:r-admm-five-eqs}
The implementation of the R-ADMM algorithm described in the set of five equations given in \eqref{eq:psi-update}, \eqref{eq:xi-update} and \eqref{eq:prs-3} applied to the dual of problem \eqref{eq:primal-indicator-f}, reduces to alternating between the following two updates

\begin{align}\label{eq:x-update-distributed}
	& x_i(k)=\argmin_{x_i}\left\{f_i(x_i)-\left(\sum_{j\in\mathcal{N}_i}z_{ji}^\top(k)\right) x_i\right.\\
	&\,\,\,\qquad\qquad\qquad \qquad \qquad \qquad\qquad    \Biggl.+\,\frac{\rho}{2}|\mathcal{N}_i|\|x_i\|^2 \Biggr\},\nonumber
\end{align}
for all $i \in V$, and 
\begin{align}\label{eq:z-update-distributed}
\begin{split}
	& z_{ij}(k+1)=(1-\alpha)z_{ij}(k)-\alpha z_{ji}(k)+2\alpha\rho x_i(k)\\
	& z_{ji}(k+1)=(1-\alpha)z_{ji}(k)-\alpha z_{ij}(k)+2\alpha\rho x_j(k)
\end{split}
\end{align}
for all $(i,j) \in \mathcal{E}$.\oprocend
\end{prop}

\smallskip

\begin{remark}
Observe that the reformulation of \eqref{eq:psi-update}, \eqref{eq:xi-update} and \eqref{eq:prs-3} as in Proposition~\ref{pr:r-admm-five-eqs} is possible for the particular structure of Problem~\eqref{eq:primal-indicator-f} and, in particular, for the structure of the constraints $Ax+y=0$. In general, given a set of constraints $Ax+By=c$ being $A$, $B$ and $c$ generic matrices and vector, such reformulation might not be possible.\oprocend
\end{remark}

\smallskip

The previous proposition naturally suggests an alternative distributed implementation of the R-ADMM Algorithm~\ref{alg:r-admm-three-eqs}, in which each node $i$ stores in its local memory the variables $x_i$ and $z_{ij},j\in\mathcal{N}_i$. Then, at each iteration of the algorithm, each node $i$ first collects the variables $z_{ji},j\in\mathcal{N}_i$; second, updates $x_i$ and $z_{ij}$ according to \eqref{eq:x-update-distributed} and the first of \eqref{eq:z-update-distributed}, respectively; finally, it sends $z_{ij}$ to $j\in\mathcal{N}_i$.\\ 
Differently to the natural implementation just briefly described, we present a slightly different implementation building upon the observation that each node $i$, to update $x_i$ as in \eqref{eq:x-update-distributed} requires the variables $z_{ji}$ rather than $z_{ij}$ for $j\in\mathcal{N}_i$. Consequently, we assume node $i$ stores in its memory and is in charge for the update of $z_{ji},j\in\mathcal{N}_i$. The implementation is described in Algorithm~\ref{alg:smart-distributed-r-admm}.
%
%

\begin{algorithm}[ht!]
	\SetKwInOut{Input}{Input}
	\Input{Set the termination condition $K>0$. For each node $i$, initialize $x_i(0)$ and $z_{ji}(0)$, $j \in \mathcal{N}_i$.}
	$k\leftarrow0$\;
	\While{$k<K$ each agent $i$}{
		compute $x_i(k)$ according to \eqref{eq:x-update-distributed}\;
		for all $j \in \mathcal{N}_i$, compute the quantity $q_{i \to j}$ as 
		\begin{equation}\label{eq:q}
	         q_{i \to j}=-z_{ji}(k)+2 \rho x_i(k)
		\end{equation}
		for all $j \in \mathcal{N}_i$, transmit $q_{i \to j}$ to node $j$\;
		gather $q_{j \to i}$ from each neighbor $j$\;
		update $z_{ji}$ as 
		\begin{equation}\label{eq:Update_z_ji}
		z_{ji}(k+1)=(1-\alpha)z_{ji}(k)+ \alpha q_{j \to i};
		\end{equation}
		$k\leftarrow k+1$\;
	}
\caption{Modified distributed R-ADMM.}
\label{alg:smart-distributed-r-admm}
\end{algorithm}

\smallskip

As we can see, both Algorithms~\ref{alg:r-admm-three-eqs} and \ref{alg:smart-distributed-r-admm} need a single round of transmissions at each time $k$. However, they differ for the number of packets that each node has to transmit  and for the number of variables that a node has to update. Table \ref{tab:variables-counts} reports the comparison between the two algorithms.
\begin{table}[b!]
\caption{Comparison of R-ADMM implementations.}
\label{tab:variables-counts}
\centering
\normalsize
\begin{tabular}{|c|c||c|}
\hline
& Alg. \ref{alg:r-admm-three-eqs} & Alg. \ref{alg:smart-distributed-r-admm}\\
\hline
Update and Send & $2|\mathcal{N}_i|+1$ & $|\mathcal{N}_i|+1$\\
\hline
Store & $3|\mathcal{N}_i|$ & $|\mathcal{N}_i|$ \\
\hline
\end{tabular}
\end{table}
Therefore, exploiting the auxiliary $z$ variables we have obtained an algorithm with smaller memory and computational requirements.\\
%
%
%
We conclude this section by stating the convergence properties of Algorithms~\ref{alg:r-admm-three-eqs}, \ref{alg:smart-distributed-r-admm}. The proof can be found in Appendix~\ref{app:proof-convergence}.

\smallskip

\begin{prop}\label{prop:convergence}
Consider Algorithm \ref{alg:smart-distributed-r-admm}. Let $(\alpha, \rho)$ be such that $0<\alpha <1$ and $\rho >0$. Then, for any initial conditions, the trajectories $k \to x_i(k)$, $i \in V$, generated by Algorithm~\ref{alg:smart-distributed-r-admm}, converge to the optimal solution of \eqref{eq:opt_problem}, i.e.,
$$
\lim_{k \to \infty} x_i(k) = x^*, \qquad \forall i \in \mathcal{V},
$$
for any $x_i(0)$ and $z_{ji}(0)$, $j \in \mathcal{N}_i$. The same result holds true also for  Algorithms \ref{alg:r-admm-three-eqs}.\oprocend
\end{prop}

\section{Distributed R-ADMM over lossy networks}\label{sec:robustADMM}
The distributed algorithms illustrated in the previous section work under the standing assumption that the communication channels are reliable, that is, no packet losses occur.
The goal of this section is to relax this communication requirement and, in particular, to show that Algorithm \ref{alg:smart-distributed-r-admm} still converges, under a probabilistic assumption on communication failures which is next stated. 

\smallskip

\begin{assumption}\label{ass:lossy}
During any iteration of Algorithm \ref{alg:smart-distributed-r-admm}, the communication from node $i$ to node $j$ can be lost with some probability $p$.\oprocend
\end{assumption}

\smallskip

In order to describe the communication failure more precisely, we introduce the family of independent binary random variables $L_{ij}(k)$, $k=0,1,2,\ldots$, $i \in \mathcal{V}$, $j \in \mathcal{N}_i$, such that\footnote{We highlight that the results of this section can be extended to the case where the loss probability is different for edge.}
$$
\mathbb{P}\left[L_{ij}=1\right]=p, \qquad \mathbb{P}\left[L_{ij}=0\right]=1-p.
$$  
We emphasize the fact that independence is assumed among all $L_{ij}(k)$ as $i, j$ and $k$ vary. If the packet transmitted, during the $k$-th iteration by node $i$ to node $j$ is lost, then $L_{ij}(k)=1$, otherwise $L_{ij}(k)=0$.
\\
In this lossy scenario, Algorithm \ref{alg:smart-distributed-r-admm} is modified as shown in Algorithm~\ref{alg:robust-smart-distributed-r-admm}. 

\smallskip

\begin{algorithm}[ht!]
	\SetKwInOut{Input}{Input}
	\Input{Set the termination condition $K>0$. For each node $i$, initialize $x_i(0)$ and $z_{ji}(0)$, $j \in \mathcal{N}_i$.}
	$k\leftarrow0$\;
	\While{$k<K$ each agent $i$}{
		compute $x_i(k)$ according to \eqref{eq:x-update-distributed}\;
		for all $j \in \mathcal{N}_i$, compute the quantity $q_{i \to j}$ as 
		\begin{equation}
	         q_{i \to j}=-z_{ji}(k)+2 \rho x_i(k)
		\end{equation}
		for all $j \in \mathcal{N}_i$, transmit $q_{i \to j}$ to node $j$\;
		\uIf{for $j\in\mathcal{N}_i$, $q_{j\to i}$ is received}{
			update $z_{ji}$ as 
			\begin{equation}
			z_{ji}(k+1)=(1-\alpha)z_{ji}(k)+ \alpha q_{j \to i};
			\end{equation}
		}
		\Else{
			leave $z_{ji}$ unchanged, i.e.,
			\begin{equation}
			z_{ji}(k+1)= z_{ji}(k)
			\end{equation}
		}
		$k\leftarrow k+1$\;
	}
\caption{Robust distributed R-ADMM.}
\label{alg:robust-smart-distributed-r-admm}
\end{algorithm}

\smallskip

\noindent
In this case, at $k$-th iteration node $i$ updates $x_i$ as in \eqref{eq:x-update-distributed}. Then, for $j \in \mathcal{N}_i$, it computes $q_{i \to j}$ as in \eqref{eq:q} and transmits it to node $j$. If node $j$ receives $q_{i \to j}$, then it updates $z_{ij}$ as $z_{ij}(k+1)=(1-\alpha)z_{ij}(k)+ \alpha q_{i \to j}$, otherwise $z_{ij}$ remains unchanged, i.e., $z_{ij}(k+1)=z_{ij}(k)$.
This last step can be compactly describes as
\begin{align*}
z_{ij}(k+1)&=L_{ij}(k)z_{ij}(k) + \\
&\qquad +\left(1-L_{ij}(k)\right) \,\left( (1-\alpha)z_{ij}(k)+ \alpha q_{i \to j}\right)
\end{align*}
We have the following Proposition, whose proof is reported in Appendix~\ref{app:convergence-lossy}.

\smallskip

\begin{prop}\label{prop:convergence_lossy}
Consider Algorithm \ref{alg:robust-smart-distributed-r-admm} working under the scenario described in Assumption \ref{ass:lossy}.
Let $(\alpha, \rho)$ be such that $0<\alpha <1$ and $\rho >0$. Then, for any initial conditions, the trajectories $k \to x_i(k)$, $i \in \mathcal{V}$, generated by Algorithm~\ref{alg:robust-smart-distributed-r-admm}, converge almost surely to the optimal solution of \eqref{eq:opt_problem}, i.e.,
$$
\lim_{k \to \infty} x_i(k) = x^*, \qquad \forall i \in \mathcal{V},
$$
with probability one, for all $i \in \mathcal{V}$, for any $x_i(0)$ and $z_{ji}(0)$, $j \in \mathcal{N}_i$.\oprocend
\end{prop}

\smallskip

We stress that the underling idea behind the result of Proposition~\ref{prop:convergence_lossy} relies on rewriting Algorithm~\ref{alg:robust-smart-distributed-r-admm} as a stochastic KM iteration and then to resort to a different set of methodological tools from probabilistic analysis \cite{bianchi2016coordinate}.

\smallskip

\begin{remark}
We have restricted the analysis to the case of synchronous communication since we were mainly interested in investigating the algorithm performance in the presence of packet losses. The practically more appealing asynchronous scenario will be the focus of future research.\oprocend
\end{remark}

\smallskip

\begin{remark}
Interestingly, while the robustness result that we provide in the lossy scenario holds true for Algorithm~\ref{alg:robust-smart-distributed-r-admm}, we cannot prove the same for Algorithms~\ref{alg:r-admm-three-eqs} which, in the case of synchronous and reliable communications, is instead characterized by the same convergent behavior despite of the different communication and memory requirements.\oprocend
\end{remark}

\smallskip

\begin{remark}
Observe that Proposition~\ref{prop:convergence}, for the case of reliable communications, and Proposition~\ref{prop:convergence_lossy}, regarding the lossy scenario, share exactly the same region of convergence in the space of the parameters. This means that Algorithm~\ref{alg:smart-distributed-r-admm} remains provably convergent if $0<\alpha<1$ and $\rho>0$ in both cases. However, observe that the result is not \emph{necessary and sufficient} and, in particular, the convergence might hold also for value of $\alpha\geq 1$. Indeed, in the simulation Section~\ref{sec:simulation} we show that, for the case of quadratic functions $f_i,\ i\in\mathcal{V}$, the region of attraction in parameter space is larger. Moreover, despite what suggested by the intuition, the larger the packet loss probability $p$, the larger the region of convergence. However, this increased region of stability is counterbalanced by a slower convergence rate of the algorithm.\oprocend
\end{remark}

\section{Simulations}\label{sec:simulation}
\begin{figure}[t]
      \centering
      \includegraphics[width=\columnwidth]{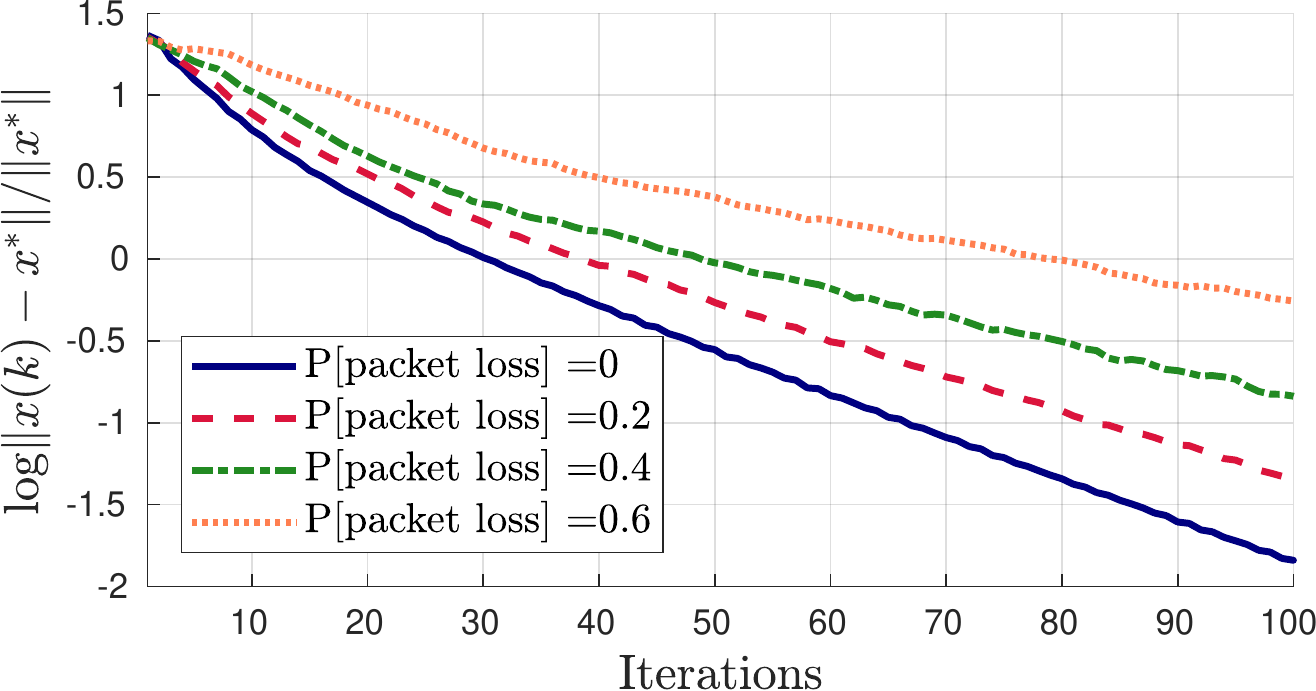}
      \caption{Evolution, in log-scale, of the relative error of Alg.~\ref{alg:robust-smart-distributed-r-admm} computed w.r.t. the unique optimal solution $x^*$ as function of different values of packet loss probability $p$ for step size $\alpha=1$ and penalty $\rho=1$. Average over 100 Monte Carlo runs.}
      \label{fig:evolution_different_losses}
   \end{figure}

\begin{figure}[t]
      \centering
      \includegraphics[width=\columnwidth]{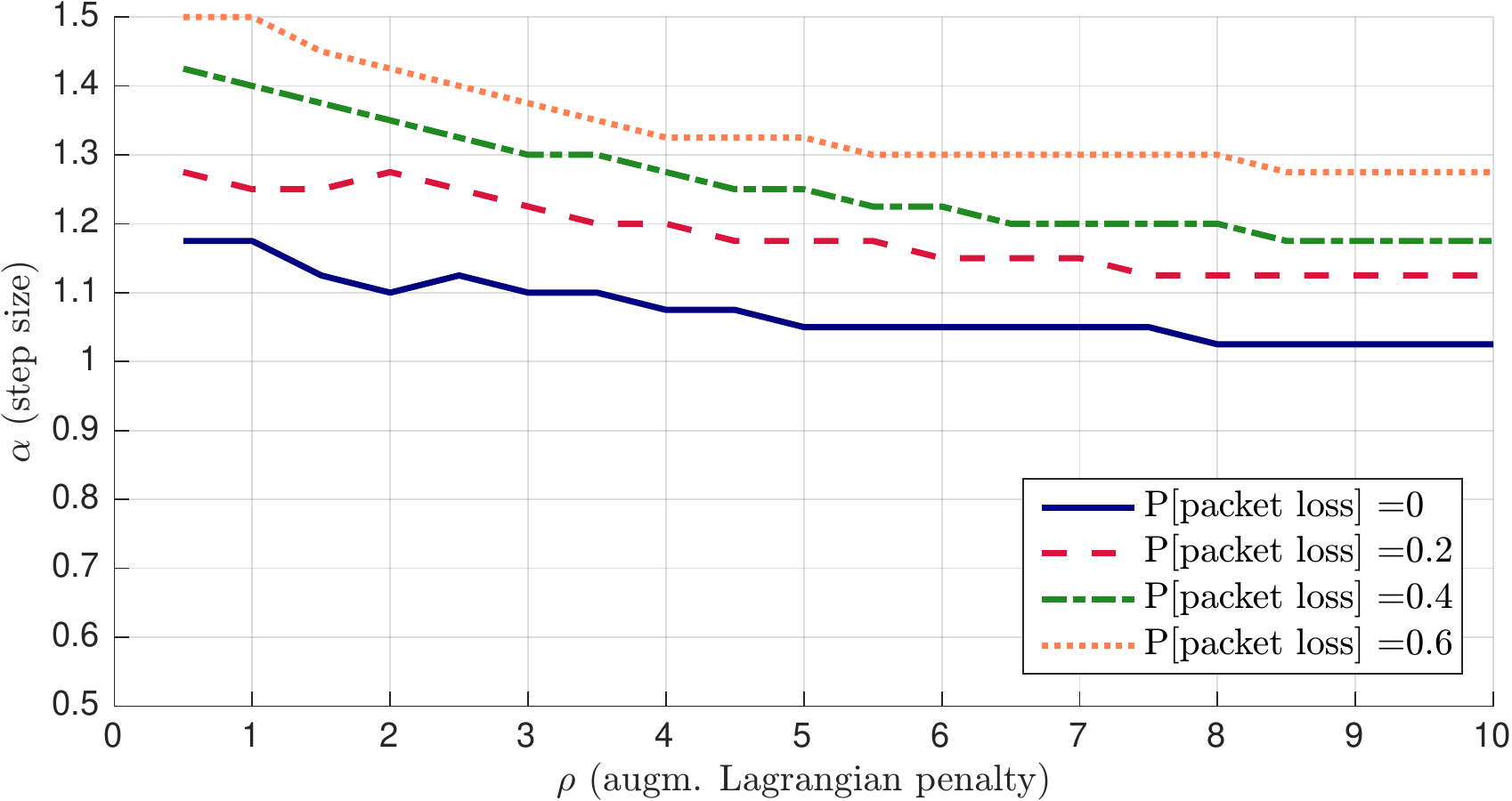}
      \caption{Stability boundaries of Alg.~\ref{alg:robust-smart-distributed-r-admm} as function of the step size $\alpha$ and the penalty $\rho$ for different values of loss probability $p$ for the family of random geometric graphs. Average over 100 Monte Carlo runs.}
      \label{fig:randgeom_stability_boundaries}
   \end{figure}

In this section we provide some experimental simulations to test the proposed R-ADMM Algorithm~\ref{alg:robust-smart-distributed-r-admm} to solve distributed consensus optimization problems \eqref{eq:opt_problem}. We are particularly interested in showing the algorithm performances in the presence of packet losses in the communication among neighboring nodes. To simplify the numerical analysis we restrict to the case of quadratic cost functions of the form
\begin{equation*}
	f_i(x_i)=a_ix_i^2 + b_ix_i + c_i
\end{equation*}
where, in general, the quantities $a_i,b_i,c_i\in\mathbb{R}$ are different for each node $i$. In this case the update of the primal variables becomes linear and, in particular, Eq.~\eqref{eq:x-update-distributed} reduces to
\begin{align*}
	&x_i(k)=\frac{\sum_{j\in\mathcal{N}_i}z_{ji}(k)-b_i}{2a_i+\rho|\mathcal{N}_i|}\, .
\end{align*}
%
We consider the family of random geometric graphs with $N=10$ and communication radius $r=0.1$[p.u.] in which two nodes are connected if and only if their relative distance is less that $r$.
We perform a set of 100 Monte Carlo runs for different values of packet losses probability $p$, step size $\alpha$ and penalty parameters $\rho$.\\
First of all, for different values of packet loss probability $p$ and for fixed values of step size $\alpha=1$ and penalty $\rho=1$, Figure~\ref{fig:evolution_different_losses} shows the evolution of the relative error
$$
\log\frac{\|x(k)-x^*\|}{\|x^*\|}
$$
computed with respect to the unique minimizer $x^*$ and averaged over 100 Monte Carlo runs. As expected, the higher the packet loss probability, the smaller the rate of convergence. Indeed, failures in the communication among neighboring nodes negatively affect the computations.\\
Figure~\ref{fig:randgeom_stability_boundaries} plots the stability boundaries of the R-ADMM Algorithm~\ref{alg:robust-smart-distributed-r-admm} as function of step size $\alpha$ and penalty $\rho$ for different packet loss probabilities $p$. More specifically, each curve in Figure~\ref{fig:randgeom_stability_boundaries} represents the numerical boundary below which the algorithm is found to be convergent and above which, conversely, the algorithm diverges. In this case the results turn out extremely interesting. Indeed, given $\alpha$ and $\rho$, for increasing packet loss probability $p$, the stability region enlarges. This  means that the higher the loss probability is, the more robust the algorithm is. The numerical findings are perfectly in line with the result of Proposition~\ref{prop:convergence_lossy}, telling us that for $\alpha\in (0,1)$ the algorithm converges for any value of $\rho$. However, it suggests the additional interesting fact that the theory misses to capture a larger area -- in parameters space and depending on $p$ -- for which the algorithm still converges. This will certainly be a direction of future  investigation.

\begin{figure}[t]
      \centering
      \includegraphics[width=\columnwidth]{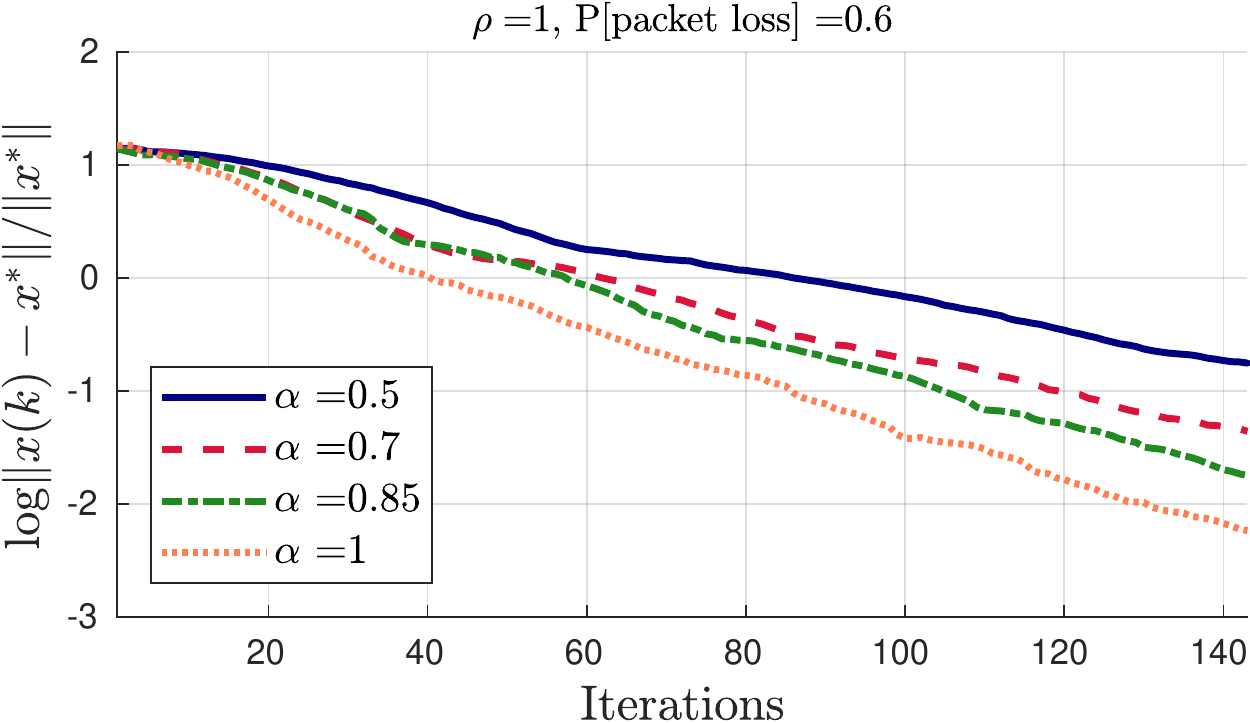}
      \caption{Evolution, in log-scale, of the relative error of Alg.~\ref{alg:robust-smart-distributed-r-admm} computed w.r.t. the unique optimal solution $x^*$ as function of different values of the step size $\alpha$, with fixed packet loss probability $p=0.6$ and penalty $\rho=1$. Average over 100 Monte Carlo runs.}
      \label{fig:randgeom_different_stepsizes}
   \end{figure}

Finally, Figure~\ref{fig:randgeom_different_stepsizes} reports the evolution of the error as a function of different values of the step-size $\alpha$. Notice that to values of $\alpha$ that are larger than $1/2$ correspond faster convergences. Recalling that setting $\alpha=1/2$ yields the standard ADMM, then it is clear that the use of the R-ADMM can speed up the convergence, which motivates its use against the use of the classic ADMM.

\section{Conclusions and Future Directions}\label{sec:conclusions}
In this paper we addressed the problem of distributed consensus optimization in the presence of synchronous but unreliable communications. 
Building upon results in operator theory on Hilbert spaces, we leveraged the relaxed Peaceman-Rachford Splitting operator to introduced what is referred to R-ADMM, a generalization of the well known ADMM algorithm. We started by drawing some interesting connections with the classical formulation as typically presented. Then, we introduced several algorithmic reformulations of the R-ADMM which differs in terms of computational, memory and communication requirements. Interestingly the last implementation, besides being extremely light from both the communication and memory point of views, turns out the be provably robust to random communication failures. Indeed, we rigorously proved how, in the lossy scenario, the region of convergence in parameters space remains unchanged compared to the case of reliable communication; yet, we numerically showed that the region of convergence is positively affected by a larger packet loss probability. The drawback lies in a slower convergence rate of the algorithm.\\
There remain many open questions paving the paths to future research directions such as analysis of the asynchronous case and generalization of the results to more general distributed optimization problems.






\appendices

\section{Derivation of Algorithm \ref{alg:r-admm-three-eqs}}\label{app:derivation-alg-3eqs}
First of all we derive the augmented Lagrangian \eqref{eq:augmented-lagr} for problem \eqref{eq:primal-indicator-f}, and obtain
\begin{align}\label{eq:augmented-lagr-distributed}
\begin{split}
	\mathcal{L}_\rho(x,y;w)=\sum_{i=1}^Nf_i(x_i)&+\iota_{(I-P)}(y)+\\&-w^\top(Ax+y)+\frac{\rho}{2}\|Ax+y\|^2,
\end{split}
\end{align}
where $\|Ax+y\|^2=\|Ax\|^2+\|y\|^2+2\langle Ax,y\rangle$. We can now proceed to derive equations \eqref{eq:r-admm-1}--\eqref{eq:r-admm-3} for the problem at hand.

\subsubsection{Equation \eqref{eq:r-admm-1}}
By \eqref{eq:augmented-lagr-distributed} and discarding the terms that do not depend on $y$ we get
\begin{align*}
	y(k+1)=\argmin_y&\Big\{\iota_{(I-P)}(y)-w^\top(k)y+\frac{\rho}{2}\|y\|^2\\&+2\alpha\rho\langle Ax(k),y\rangle+\rho(2\alpha-1)\langle y,y(k)\rangle\Big\}
\end{align*}
where we summed the terms with the inner product $\langle Ax(k),y\rangle$. Therefore we need to solve the problem
\begin{align*}
	y(k+1)=\argmin_{y=Py}&\Big\{-w^\top(k)y+\frac{\rho}{2}\|y\|^2 \\
	&+2\alpha\rho\langle Ax(k),y\rangle+\rho(2\alpha-1)\langle y,y(k)\rangle\Big\}
\end{align*}
that for simplicity we can write as
\begin{equation}\label{eq:problem-y}
	y(k+1)=\argmin_{y=Py}\{h_{\alpha,\rho}(y;x(k),w(k))\}.
\end{equation}\\
We apply now the Karush-Kuhn-Tucker (KKT) conditions \cite{boyd2004convex} to problem \eqref{eq:problem-y} and obtain the system
\begin{align}
	&\nabla\Big[h_{\alpha,\rho}(y;x(k),w(k))-\nu^\top(I-P)y\Big|_{y(k+1),\nu^*}=0\label{eq:kkt-1}\\
	&y(k+1)=Py(k+1)\label{eq:kkt-2}
\end{align}
where $\nu^*$ is the optimal value of the Lagrange multipliers of the problem.\\
By computing the gradient in \eqref{eq:kkt-1} we obtain
\begin{align}\label{eq:kkt-1-bis}
\begin{split}
	y(k+1)=\frac{1}{\rho}\big[w(k)&-2\alpha\rho Ax(k)\\&-\rho(2\alpha-1)y(k)+(I-P)\nu^*\big].
\end{split}
\end{align}
We substitute this formula for $y(k+1)$ in the right-hand side of \eqref{eq:kkt-2} which results in
\begin{align}\label{eq:kkt-2-bis}
\begin{split}
	y(k+1)=\frac{1}{\rho}\big[P&w(k)-2\alpha\rho PAx(k)\\&-\rho(2\alpha-1)Py(k)-(I-P)\nu^*\big]
\end{split}
\end{align}
for the fact that $P^2=I$ and hence $P(I-P)=-(I-P)$.\\
We sum now equations \eqref{eq:kkt-1-bis} and \eqref{eq:kkt-2-bis} and obtain
\begin{align}\label{eq:y-update-final}
\begin{split}
	y(k+1)=\frac{1}{2\rho}(I+P)\big[w(k)-2\alpha\rho Ax(k)-\rho(2\alpha-1)y(k)\big].
\end{split}
\end{align}
Finally noting that, given a vector $t$ of dimension equal to that of $y$, the $ij$-th element of $(I+P)t$ is equal to $t_{ij}+t_{ji}$, then the update for $y_{ij}(k+1)$ follows.

\subsubsection{Equation \eqref{eq:r-admm-2}}
By equation \eqref{eq:r-admm-2} and \eqref{eq:y-update-final} we can write
\begin{align*}
	w(k+1)=&w(k)-2\alpha\rho Ax(k)-\rho(2\alpha-1)y(k)+\\
	&-\frac{1}{2}(I+P)[w(k)-2\alpha\rho Ax(k)-\rho(2\alpha-1)y(k)]\\
	=&\frac{1}{2}(I-P)[w(k)-2\alpha\rho Ax(k)-\rho(2\alpha-1)y(k)]
\end{align*}
and by the definition of $I-P$ we get the update equation for $w_{ij}(k+1)$ stated in Algorithm \ref{alg:r-admm-three-eqs}.

\subsubsection{Equation \eqref{eq:r-admm-3}}
Finally we apply equation \eqref{eq:r-admm-3} to the problem at hand, which means that we need to solve
\begin{align*}
	x(k+1)=&\argmin_x\Bigg\{\sum_{i=1}^Nf_i(x_i)+\\&-\Big(w(k+1)-\rho y(k+1)\Big)^\top Ax+\frac{\rho}{2}\|Ax\|^2\Bigg\}.
\end{align*}
We know that each variable $x_i$ appears in $|\mathcal{N}_i|$ constraints and therefore $\|Ax\|^2=\sum_{i=1}^N|\mathcal{N}_i|\|x_i\|^2$. Moreover, given a vector $t$ with the same size as $y$, we have
\begin{align*}
	t^\top Ax&=
		\begin{bmatrix}
			\cdots & t_{ji}^\top & \cdots & t_{ji}^\top & \cdots
		\end{bmatrix}
		\begin{bmatrix}
			\vdots\\-x_i\\ \vdots\\-x_j\\ \vdots
		\end{bmatrix}\\
	&=\sum_{(i,j)\in\mathcal{E}}\left(t_{ji}^\top x_i+t_{ij}^\top x_j\right)\\
	&=\sum_{i=1}^N\left(\sum_{j\in\mathcal{N}_i}t_{ji}^\top\right)x_i.
\end{align*}
and we get the update equation for $x_i(k+1)$ substituting $\Big(w(k+1)-\rho y(k+1)\Big)$ to $t$.
Notice that by the results obtained above we have
\begin{align*}
	\Big(w(k+1)-&\rho y(k+1)\Big)=\\&=-P[w(k)-2\alpha\rho Ax(k)-\rho(2\alpha-1)y(k)]
\end{align*}
which means that $x(k+1)$ can be computed as a function of the $x$, $y$ and $w$ variables at time $k$ only.\oprocendbis

\section{Proof of Proposition \ref{pr:r-admm-five-eqs}}\label{app:proof-prop-1}

\subsubsection{Equations \eqref{eq:psi-update}}
The following derivation shares some points with the derivation described in the section above. Indeed, applying the first equation of \eqref{eq:psi-update} to the problem at hand requires that we solve
\begin{equation*}
	y(k)=\argmin_{y=Py}\left\{-z^\top(k)y+\frac{\rho}{2}\|y\|^2\right\},
\end{equation*}
which can be done by solving the system of KKT conditions of the problem as performed above. The result is
\begin{equation}\label{eq:y-update-proof}
	y(k)=\frac{1}{2\rho}(I+P)z(k).
\end{equation}
It easily follows from \eqref{eq:y-update-proof} that $\psi(k)=\frac{1}{2}(I-P)z(k)$.

\subsubsection{Equations \eqref{eq:xi-update}}
First of all we have $(2\psi(k)-z(k))=-Pz(k)$, hence according to the same reasoning employed above to derive the expression for $x(k+1)$ we find \eqref{eq:x-update-distributed}. Moreover, we have $\xi(k)=-Pz(k)-\rho Ax(k)$.

\subsubsection{Equation \eqref{eq:prs-3}}
By the results derived above we can easily compute
\begin{equation*}
	z(k+1)=(1-\alpha)z(k)-\alpha Pz(k)-2\alpha\rho Ax(k)
\end{equation*}
which gives equations \eqref{eq:z-update-distributed}.

\smallskip

Notice that to compute the variables $y(k)$, $\psi(k)$, $x(k)$ and $\xi(k)$ we need only the variables $z(k)$. Moreover, to update $z$ we require only $z(k)$ and $x(k)$. Hence the five update equations reduce to the updates for $x$ and $z$ only.\oprocendbis

\section{Proof of Proposition \ref{prop:convergence}}\label{app:proof-convergence}
To prove convergence of the R-ADMM in the two implementations of Algorithms \ref{alg:r-admm-three-eqs} and \ref{alg:smart-distributed-r-admm}, we resort to the following result, adapted from \cite[Corollary 27.4]{bauschke2011convex}.
\begin{prop}[{\cite[Corollary 27.4]{bauschke2011convex}}]\label{pr:convergence-deterministic}
Consider problem \eqref{eq:MinProblem} and assume that it has solution; let $\alpha\in(0,1)$, $\rho>0$, and $x(0)\in\mathcal{X}$. Assume to apply equations \eqref{eq:prs-1}--\eqref{eq:prs-3} to the problem. Then there exists $z^*$ such that
\begin{itemize}
	\item $x^*=\prox_{\rho g}(z^*)\in\argmin_x\{f(x)+g(x)\}$, and
	\item $\{z(k)\}_{k\in\mathbb{N}}$ converges weakly to $z^*$.
\end{itemize}\oprocend
\end{prop}

We need to show now that this result applies to the dual problem of problem \eqref{eq:primal-indicator-f}. First of all, by formulation of the problem we have that $f$ is convex and proper (and also closed). Moreover, by \cite[Example 8.3]{bauschke2011convex} we know that the indicator function of a convex set is convex (and, by definition, proper). But the set of vectors $y$ that satisfy $(I-P)y=0$ is indeed convex, hence also $g$ is convex and proper.\\
Now \cite[Theorem 12.2]{rockafellar2015convex} states that the convex conjugate of a convex and proper function is closed, convex and proper. Therefore both $d_f$ and $d_g$ are closed, convex and proper, which means that we can apply the convergence result in Proposition \ref{pr:convergence-deterministic} to the dual problem of \eqref{eq:primal-indicator-f}.\\
Therefore we have that $w^*=\prox_{\rho d_g}(z^*)$ is indeed a solution of the dual problem and $\{z(k)\}_{k\in\mathbb{N}}$ converges to $z^*$. But since the duality gap is zero, then when we attain the optimum of the dual problem we have obtained that of the primal as well.\oprocendbis

\section{Proof of Proposition \ref{prop:convergence_lossy}}\label{app:convergence-lossy}
In order to prove the convergence of Algorithm \ref{alg:robust-smart-distributed-r-admm} we need to introduce a probabilistic framework in which to reformulate the KM update. For this stochastic version of the KM iteration we can state a convergence result adapted from \cite[Theorem 3]{bianchi2016coordinate} and show that indeed Algorithm \ref{alg:robust-smart-distributed-r-admm} is represented by this formulation.
 
We are therefore interested in altering the standard KM iteration \eqref{eq:km-iteration} in order to include a stochastic selection of which coordinates in $\mathcal{I}=\{1,\ldots,M\}$ to update at each instant. To do so we introduce the operator $\hat{T}^{(\xi)}:\mathcal{X}\rightarrow\mathcal{X}$ whose $i$-th coordinate is given by $\hat{T}^{(\xi)}_ix=T_ix$ if the coordinate is to be updated ($i\in \xi$), $\hat{T}^{(\xi)}_ix=x_i$ otherwise ($i\not\in \xi$). In general the subset of coordinates to be updated changes from one instant to the next. Therefore, on a probability space $(\Omega,\mathcal{F},\mathbb{P})$, we define the random i.i.d. sequence $\{\xi_k\}_{k\in\mathbb{N}}$, with $\xi_k:\Omega\rightarrow 2^\mathcal{I}$, to keep track of which coordinates are updated at each instant.
The stochastic KM iteration is finally defined as
\begin{equation}\label{eq:stochastic-km}
	x(k+1)=(1-\alpha)x(k)+\alpha\hat{T}^{(\xi_{k+1})}x(k)
\end{equation}
and consists of the $\alpha$-averaging of a stochastic operator.

The stochastic iteration satisfies the following convergence result, which is particularized from \cite{bianchi2016coordinate} using the fact that a nonexpansive operator is $1$-averaged, and a constant step size.
\begin{prop}[{\cite[Theorem 3]{bianchi2016coordinate}}]
Let $T$ be a nonexpansive operator with at least a fixed point, and let the step size be $\alpha\in(0,1)$. Let $\{\xi_k\}_{k\in\mathbb{N}}$ be a random i.i.d. sequence on $2^\mathcal{I}$ such that
\begin{equation*}
	\forall i\in\mathcal{I},\ \exists I\in2^\mathcal{I}\ \text{s.t.}\ i\in I\ \text{and}\ \mathbb{P}[\xi_1=I]>0.
\end{equation*}
Then for any deterministic initial condition $x(0)$ the stochastic KM iteration \eqref{eq:stochastic-km} converges almost surely to a random variable with support in the set of fixed points of $T$.\oprocend
\end{prop}

We turn now to the distributed optimization problem, in which the stochastic KM iteration is performed on the auxiliary variables $z$. In particular we assume that the packet loss occurs with probability $p$, and that in the case of packet loss the relative variable is not updated. As shown in the main paper, this update rule can be compactly written as
\begin{equation}\label{eq:operator-packet-loss}
	\hat{T}^{(\xi_{k+1})}z(k)=L_kz(k)+(I-L_k)Tz(k)
\end{equation}
where $L_k$ is the diagonal matrix with elements the realizations of the binary random variables that model the packet loss at time $k$. Recall that these variables take value 1 if the packet is lost.\\
Substituting now the operator \eqref{eq:operator-packet-loss} into \eqref{eq:stochastic-km} we get the update equation
\begin{equation}\label{eq:stochastic-km-order-1}
	z(k+1)=(1-\alpha)z(k)+\alpha\left[L_kz(k)+(I-L_k)Tz(k)\right]
\end{equation}
which conforms to the stochastic KM iteration for which the convergence result is stated.\\
Finally, notice that in the main article the $\alpha$-averaging is applied before the stochastic coordinate selection, that is the update is given by
\begin{equation}\label{eq:stochastic-km-order-2}
	z(k+1)=L_kz(k)+(I-L_k)\left[(1-\alpha)z(k)+\alpha Tz(k)\right].
\end{equation}
However it can be easily shown that \eqref{eq:stochastic-km-order-1} and \eqref{eq:stochastic-km-order-2} do indeed coincide, hence proving the convergence of our update scheme.\oprocendbis

\bibliographystyle{./IEEEtran} 
\bibliography{./IEEEabrv,./references}

\end{document}